\documentclass[12pt,twoside]{article}

\usepackage{booktabs,subcaption,amsfonts,dcolumn}
\newcolumntype{d}[1]{D..{#1}}
\newcommand\mc[1]{\multicolumn{1}{c}{#1}} 

\usepackage{amsmath,amssymb,amsthm,amscd,color,accents,mathrsfs,url}

\usepackage{enumitem}
\usepackage{euscript}
\usepackage{multirow}
\usepackage{latexsym,amsbsy}
\usepackage[pdftex]{graphicx}
\usepackage{epsfig}
\usepackage{subdepth}

\usepackage{float}
\usepackage{color}
\usepackage{relsize}
\usepackage[geometry]{ifsym}
\usepackage{pdfpages}
\usepackage[makeroom]{cancel}

\usepackage[pagebackref=false,colorlinks=true,linkcolor=blue,citecolor=blue,urlcolor=blue]{hyperref}

\numberwithin{equation}{section}

\everymath{\displaystyle}

\theoremstyle{plain}
\newtheorem{theorem}{Theorem}[section]
\newtheorem{proposition}[theorem]{Proposition}
\newtheorem{lemma}[theorem]{Lemma}
\newtheorem{corollary}[theorem]{Corollary}

\makeatletter
\newcommand*\bigcdot{\mathpalette\bigcdot@{.8}}
\newcommand*\bigcdot@[2]{\mathbin{\vcenter{\hbox{\scalebox{#2}{{\hskip 1pt}$\m@th#1\bullet$}}}}}
\makeatother

\def\E {\mathbb{E}}

\def\R {\mathbb{R}}

\def\Cov {\mathrm{Cov}}

\def\dd {\hskip2pt\mathrm{d}}
\def\tr {\hskip1pt\mathrm{tr}\,}

\def\mcalS{\mathcal{S}}

\newcommand{\Proof}{\noindent\textit{Proof.}\ \ }

\begin{document}

\title{
\Large\textbf{Complete Asymptotic Expansions and the High-Dimensional Bingham Distributions}}

\author{{Armine Bagyan}\thanks{Department of Statistics, Pennsylvania State University, University Park, PA 16802, U.S.A. \ E-mail address: aub171@psu.edu}
\ {and Donald Richards}\thanks{Department of Statistics, Pennsylvania State University, University Park, PA 16802, U.S.A. \ E-mail address: richards@stat.psu.edu
\endgraf
\ $^\dag$Corresponding author.
}}

\date{\today}

\maketitle

\begin{abstract}
For $d \ge 2$, let $X$ be a random vector having a Bingham distribution on $\mcalS^{d-1}$, the unit sphere centered at the origin in $\R^d$, and let $\Sigma$ denote the symmetric matrix parameter of the distribution.  Let $\Psi(\Sigma)$ be the normalizing constant of the distribution and let $\nabla \Psi_d(\Sigma)$ be the matrix of first-order partial derivatives of $\Psi(\Sigma)$ with respect to the entries of $\Sigma$.  We derive complete asymptotic expansions for $\Psi(\Sigma)$ and $\nabla \Psi_d(\Sigma)$, as $d \to \infty$; these expansions are obtained subject to the growth condition that $\|\Sigma\|$, the Frobenius norm of $\Sigma$, satisfies $\|\Sigma\| \le \gamma_0 d^{r/2}$ for all $d$, where $\gamma_0 > 0$ and $r \in [0,1)$.  Consequently, we obtain for the covariance matrix of $X$ an asymptotic expansion up to terms of arbitrary degree in $\Sigma$.  Using a range of values of $d$ that have appeared in a variety of applications of high-dimensional spherical data analysis we tabulate the bounds on the remainder terms in the expansions of $\Psi(\Sigma)$ and $\nabla \Psi_d(\Sigma)$ and we demonstrate the rapid convergence of the bounds to zero as $r$ decreases.



\smallskip
\noindent
{\em Running head}: The High-Dimensional Bingham Distribution.
\end{abstract}

\section{Introduction}
\label{sec_introduction}

Let $d$ be a positive integer, where $d \ge 2$.  Denote by $\mcalS^{d-1}$ the unit sphere in $d$-dimensional Euclidean space, $\R^d$, with center at the origin, and let $\dd x$ be the normalized surface measure on $\mcalS^{d-1}$.  For a symmetric $d \times d$ matrix $\Sigma = (\sigma_{i,j})$, the $d \times 1$ random vector $X = (X_1,\ldots,X_d)'$ is said to have a \textit{Bingham distribution} on $\mcalS^{d-1}$ if the probability density function of $X$, with respect to the surface measure $\dd x$, is 
\begin{equation}
\label{eq_Bingham_pdf}
\phi(x;{\hskip-3pt}\Sigma) = [\Psi_d(\Sigma)]^{-1} \exp(x'\Sigma x),
\end{equation}
$x \in \mcalS^{d-1}$, with normalizing constant 
\begin{equation}
\label{eq_normaliz_const}
\Psi_d(\Sigma) = \int_{\mcalS^{d-1}} \exp(x' \Sigma x) \dd x.
\end{equation}
These distributions were introduced by Bingham \cite{Bingham}, and they have become prominent since their appearance.  We cite the work of Bagyan \cite{Bagyan}, Bhattacharya and Bhattacharya \cite{Bhattacharya}, Bingham, Chang, and Richards \cite{BinghamCR}, Chikuse \cite{Chikuse}, Dryden \cite{Dryden}, Dryden and Mardia \cite{DrydenM}, Fisher, Lewis, and Embleton \cite{Fisher}, Kume and Wood \cite{KumeW}, Mardia and Jupp \cite{MardiaJ}, and Patrangenaru and Ellingson \cite{Patrangenaru} for various aspects of the theory of these distributions, their applications in random modulation, spherical data analysis, Procrustes analysis, shape analysis, and for many references to the literature.  

One approach in the study of high-dimensional approximations for the Bingham distribution is to approximate the distribution of $X$ and then to deduce approximations for, e.g., $\Cov(X)$, the covariance matrix of $X$; see, e.g., Chikuse \cite{Chikuse}, Dryden \cite{Dryden}, and Kume and Wood \cite{KumeW}.  In this article, starting from the zonal polynomial series for $\Psi_d(\Sigma)$, we obtain complete asymptotic expansions, i.e., series expansions in negative powers of $d$, for $\Psi_d(\Sigma)$ and $\Cov(X)$ as $d \to \infty$.  An advantage of our approach is that the resulting asymptotic expansions are obtained purely in terms of $\|\Sigma\| = [\tr(\Sigma^2)]^{1/2}$, the Frobenius norm of $\Sigma$.  Thus we make no assumptions about the structure of individual matrix entries of $\Sigma$; indeed, our only constraint on $\Sigma$ is that $\|\Sigma\| \le \gamma_0 d^{r/2}$ for all $d$, with universal constants $\gamma_0 > 0$ and $r \in [0,1)$.   

A key role is played in our derivations by masterful results of Reznick \cite{Reznick}, who obtained bounds on the ratios of power sum symmetric functions in several variables.  Reznick's results enable us to obtain, for series involving $\Sigma$, bounds which are valid for all $d$.  

In Section \ref{sec_asymp_normaliz_const}, we derive a complete asymptotic expansion for $\Psi_d(\Sigma)$.  We obtain in Proposition \ref{prop_Ck_series_tail_bound} a $O(\cdot)$ bound for the remainder of the zonal polynomial series expansion of $\Psi_d(\Sigma)$ upon truncation after any given number of terms.  

In Section \ref{sec_asymp_cov_matrix}, we obtain a complete asymptotic expansion for the matrix 
$$
\Psi_d(\Sigma) \Cov(X) = \int_{\mcalS^{d-1}} xx' \exp(x' \Sigma x) \dd x.
$$
We derive in Proposition \ref{prop_termwise_diff} a $O(\cdot)$ bound for the Frobenius norm of the remainder of that series expansion after truncation at any number of terms.  We obtain the asymptotic expansion for $\Cov(X)$ by multiplying together the series for $[\Psi_d(\Sigma)]^{-1}$ and $\Psi_d(\Sigma) \Cov(X)$, and we provide in Theorem \ref{th_cov_matrix} the explicit expansion of $\Cov(X)$ up to terms of arbitrary degree in $\Sigma$.  We remark that our methods can also be applied to obtain complete asymptotic expansions for higher moments of $X$.  

In Section \ref{sec_implications}, we describe the range of values of $d$ that have appeared in the literature in a wide variety of applications of statistical analyses of spherical data.  We compute in Tables \ref{table_moderate_r} and \ref{table_large_r} numerous values of the bounds for the remainder terms in the asymptotic expansions of $\Psi_d(\Sigma)$ and $\nabla\Psi_d(\Sigma)$, the matrix gradient of $\Psi_d(\Sigma)$, and we elucidate some observed features of the tabulated values.  Finally we provide all proofs in the appendix Section \ref{sec_appendix}.

\section{Asymptotic expansion of the normalizing constant}
\label{sec_asymp_normaliz_const}

Denote by $\R_{sym}^{d \times d}$ the vector space of $d \times d$ real symmetric matrices.  As with any norm on a vector space, the Frobenius norm, $\|\Sigma\| = [\tr (\Sigma^2)]^{1/2}$, satisfies the subadditivity property, $\|\Sigma_1 + \Sigma_2\| \le \|\Sigma_1\| + \|\Sigma_2\|$ for all $\Sigma_1, \Sigma_2 \in \R_{sym}^{d \times d}$.  The Frobenius norm is also submultiplicative: $\|\Sigma^j\| \le \|\Sigma\|^j$ for all $\Sigma \in \R_{sym}^{d \times d}$ and $j = 1,2,3,\ldots$ \cite[p.~341]{Horn}.  

Let $g_1$ be a real-valued function of $d$ and $g_2$ be a nonnegative function of $d$.  We use the terminology $g_1(d) = O(g_2(d))$ as $d \to \infty$ to signify that there exists a constant $\gamma > 0$ such that $\vert g_1(d)\vert \le \gamma g_2(d)$ for all sufficiently large $d$.  

Consider the situation in which $\|\Sigma\| \le \gamma_0$ for all $d$.  By expanding the exponential function in \eqref{eq_normaliz_const} in a Maclaurin expansion and integrating the first two terms of that series, we find that $\Psi_d(\Sigma) \approx 1 + d^{-1} \tr (\Sigma) + O(d^{-2})$.  However, if $\|\Sigma\|$ is allowed to grow as $d$ increases then a deeper analysis is required to derive the order of approximation of the remainder term.  

Thus, our governing assumption throughout the article is that $\|\Sigma\| \le \gamma_0 d^{r/2}$ for all $d$, with constants $\gamma_0 > 0$ and $r \in [0,1)$.  

For $a \in \R$, define the \textit{rising factorial} $(a)_k = a(a+1)\cdots(a+k-1)$, $k = 0,1,2,\ldots$.  The \textit{zonal polynomial} corresponding to the partition $(k)$ may be defined as 
\begin{equation}
\label{eq_Ck_multisum}
C_{(k)}(\Sigma) = \frac{k!}{(1/2)_k} \ \sum_{i_1 + 2 i_2 + 3 i_3 + \cdots + k i_k = k} \ \prod_{j=1}^k \frac{[\tr (\Sigma^j)]^{i_j}}{i_j! \, (2j)^{i_j}},
\end{equation}
where the multiple sum is over all nonnegative integers $i_1,\ldots,i_k$ such that $i_1 + 2 i_2 + 3 i_3 + \cdots + k i_k = k$ (James \cite[p.~493, Eq.~(123)]{James}).  In particular, 
\begin{equation}
\label{eq_zonal_polys}
C_{(k)}(\Sigma) = 
\begin{cases}
1, & k=0 \\
\tr(\Sigma), & k=1 \\
\frac{1}{3}\big[\big(\tr (\Sigma)\big)^2 + 2\tr(\Sigma^2)\big], & k=2
\end{cases}.
\end{equation}

It is well-known that $\Psi_d(\Sigma)$, the normalizing constant in \eqref{eq_normaliz_const}, can be expressed in terms of a confluent hypergeometric function of matrix argument.  We refer to Bingham \cite{Bingham} for the case $d=3$ and to Muirhead \cite[p.~288, Example 7.8]{Muirhead} for the general case.  The confluent hypergeometric function of matrix argument can be expressed as an infinite series of zonal polynomials, leading to the result, 
\begin{equation}
\label{eq_c_series}
\Psi_d(\Sigma) = {}_1F_1(1/2;d/2;\Sigma) 
= \sum_{k=0}^\infty \frac{(1/2)_k}{(d/2)_k} \frac{C_{(k)}(\Sigma)}{k!};
\end{equation}
see James \cite[p.~477]{James}, Muirhead \cite[Chapter 7]{Muirhead}.  The series \eqref{eq_c_series} converges absolutely for each $\Sigma$ and uniformly on compact regions (Gross and Richards \cite[Theorem 6.3]{Gross}).  

For $m \ge 1$ we will use as an asymptotic approximation to $\Psi_d(\Sigma)$ the sum of the first $m$ terms of the series \eqref{eq_c_series}, 
\begin{equation}
\label{eq_Psi_approx}
\Psi_d(\Sigma) \approx \sum_{k=0}^{m-1} \frac{(1/2)_k}{(d/2)_k} \frac{C_{(k)}(\Sigma)}{k!}.
\end{equation}
In studying the growth rate of the remainder series, 
\begin{equation}
\label{eq_remainder_term}
R_m(\Sigma) = \sum_{k=m}^\infty \frac{(1/2)_k}{(d/2)_k} \frac{C_{(k)}(\Sigma)}{k!},
\end{equation}
we shall encounter the constants 
\begin{equation}
\label{eq_constants}
\gamma_1 = \tfrac12(1+\sqrt{3}) \simeq 1.3660, \ \gamma_2 = \gamma_0\gamma_1, \ \hbox{and} \ \gamma_3 = 2^{3/2} e^{1/2}/\gamma_1 \simeq 3.4138.
\end{equation}
Then we shall establish the following result.  

\begin{proposition}
\label{prop_Ck_series_tail_bound}
Suppose that $\|\Sigma\| \le \gamma_0 d^{r/2}$, where $\gamma_0 > 0$ and $0 \le r < 1$.  Then, for all $m \ge 1$ and all $d \ge (2\gamma_2^2)^{1/(1-r)}$, 
\begin{equation}
\label{eq_Rm_bound_4}
\vert R_m(\Sigma)\vert \le \frac{\gamma_3 \gamma_2^m}{[(m+1)!]^{1/2}} d^{-m(1-r)/2}.
\end{equation}
In particular, $R_m(\Sigma) = O\big(d^{-m(1-r)/2}\big)$ as $d \to \infty$.  
\end{proposition}

We note that a consequence of the inequality \eqref{eq_Rm_bound_4} is that $R_m(\Sigma)$ converges to zero, i.e., the accuracy of the approximation \eqref{eq_Psi_approx} increases, as $d$ increases (with $m$ held fixed) or as $m$ increases (for fixed $d$).  The upper bound \eqref{eq_Rm_bound_4} also increases as $\gamma_0$ or $r$ increases; this is consistent with the fact that the condition $\|\Sigma\| \le \gamma_0 d^{r/2}$ entails an increasingly larger matrix region as $\gamma_0$ or $r$ increases.  

As a consequence of Proposition \ref{prop_Ck_series_tail_bound}, we obtain the following asymptotic expansion for $[\Psi_d(\Sigma)]^{-1}$.  

\begin{corollary}
\label{cor_inverse_normaliz_const}
Suppose that $\|\Sigma\| \le \gamma_0 d^{r/2}$, where $\gamma_0 > 0$ and $0 \le r < 1$.  Then for all $l \ge 2$ and all $d > (6\gamma_2^2)^{1/(1-r)}$, 
\begin{equation}
\label{eq_Psi_d_inverse}
[\Psi_d(\Sigma)]^{-1} = 1 - \sum_{j=1}^{l-1} \frac{(1/2)_j}{(d/2)_j} \frac{C_{(j)}(\Sigma)}{j!} + O(d^{-(1-r)}).
\end{equation}
\end{corollary}

\section{Asymptotic expansion of the covariance matrix}
\label{sec_asymp_cov_matrix}

It is well-known that $\R_{sym}^{d \times d}$, the vector space of symmetric matrices, is an inner product space with the natural inner product $\langle H_1,H_2\rangle = \tr(H_1 H_2)$, $H_1, H_2 \in \R_{sym}^{d \times d}$.  This inner product is also related to the Frobenius norm through the identity $\langle H,H\rangle = \tr(H^2) = \|H\|^2$, $H \in \R_{sym}^{d \times d}$.  

Let $C^1(\R_{sym}^{d \times d})$ be the set of all functions $f: \R_{sym}^{d \times d} \to \R$ such that all first-order partial derivatives, $\partial f(\Sigma)/\partial\sigma_{i,j}$, are continuous, where $\Sigma = (\sigma_{i,j})$.  Let $\delta_{i,j}$ denote Kronecker's delta, i.e., $\delta_{i,j} = 1$ or $0$ according as $i=j$ or $i \neq j$, respectively; then we define $\nabla$, the \textit{gradient operator}, by 
$$
\nabla f(\Sigma) = \bigg(\tfrac12(1+\delta_{i,j})\frac{\partial}{\partial\sigma_{i,j}} f(\Sigma)\bigg).
$$
The operator $\nabla$ has appeared in various articles in multivariate analysis, e.g., Hadjicosta and Richards \cite{Hadjicosta}, Richards \cite{Richards}, Sebastiani \cite{Sebastiani}, and it arises in the theory of Taylor expansions on $\R_{sym}^{d \times d}$, as follows: For $\Sigma, H \in \R_{sym}^{d \times d}$ and $f \in C^1(\R_{sym}^{d \times d})$, there holds the first-order Taylor expansion, 
\begin{equation}
\label{eq_Taylor_expansion}
f(\Sigma+H) = f(\Sigma) + \langle H,\nabla f(\Sigma)\rangle + O(\|H\|^2), \qquad \|H\| \to 0.
\end{equation}

Turning to the covariance matrix of the Bingham-distributed random vector $X$, it follows from \eqref{eq_Bingham_pdf} that $X$ has the same distribution as $-X$; therefore, by symmetry, $\E(X) = 0$.  Hence $\Cov(X)$, the covariance matrix of $X$, equals $\E(XX')$.  By \eqref{eq_deriv_exp_trace} and a straightforward justification of the interchange of derivatives and integral, we obtain 
\begin{align*}
\Cov(X) &= [\Psi_d(\Sigma)]^{-1} \int_{\mcalS^{d-1}} xx' \exp(x' \Sigma x) \dd x \\
&= [\Psi_d(\Sigma)]^{-1} \nabla \int_{\mcalS^{d-1}} \exp\big(\tr (\Sigma xx')\big) \dd x \\
&= [\Psi_d(\Sigma)]^{-1} \nabla \Psi_d(\Sigma).
\end{align*}

Our asymptotic approximations to the matrix $\nabla \Psi_d(\Sigma)$ through term-by-term differentiation of the zonal polynomial series \eqref{eq_c_series} will be justified by means of bounds for $\|\nabla C_{(k)}(\Sigma)\|$, the Frobenius norm of the matrix $\nabla C_{(k)}(\Sigma)$.  In order to calculate bounds for each $\nabla C_{(k)}(\Sigma)$, we will need the following properties of the gradient operator.  

\begin{lemma}
\label{lem_nabla_app}
Let $H$ be a $d \times d$ symmetric matrix.  Then 
\begin{equation}
\label{eq_deriv_exp_trace}
\nabla \exp\big(\tr (\Sigma H)\big) = H \exp\big(\tr (\Sigma H)\big).
\end{equation}
Also, for $k=0,1,2,\ldots$, 
\begin{equation}
\label{eq_deriv_power_trace}
\nabla [\tr (\Sigma)]^k = k [\tr (\Sigma)]^{k-1} I_d
\end{equation}
and 
\begin{equation}
\label{eq_deriv_trace_power}
\nabla \tr(\Sigma^k) = k \Sigma^{k-1}.
\end{equation}
\end{lemma}

In the sequel, we will obtain asymptotic approximations of $\Cov(X)$ involving terms of arbitrary degree in $\Sigma$.  In deriving those results we note that, as a consequence of \eqref{eq_Ck_multisum}, \eqref{eq_deriv_power_trace}, and \eqref{eq_deriv_trace_power}, the matrix $\nabla C_{(k)}(\Sigma)$ is a homogeneous polynomial of degree $k-1$ for $k \ge 1$.  

For $m \ge 2$, we use the asymptotic approximation 
\begin{equation}
\label{eq_matrix_approx}
\nabla \Psi_d(\Sigma) \approx\sum_{k=1}^{m-1} \frac{(1/2)_k}{(d/2)_k} \frac{\nabla C_{(k)}(\Sigma)}{k!}.
\end{equation}
Note that we necessarily impose the condition $m \ge 2$ since $\nabla C_{(0)}(\Sigma) = 0$, the zero matrix.  

Following on the results in Section \ref{sec_asymp_normaliz_const}, we naturally define 
$$
\nabla R_m(\Sigma) = \sum_{k=m}^\infty \frac{(1/2)_k}{(d/2)_k} \frac{\nabla C_{(k)}(\Sigma)}{k!}.
$$
To establish the accuracy of the matrix approximation \eqref{eq_matrix_approx}, we shall derive the convergence properties as $d \to \infty$ (and also as $m \to \infty$) of $\nabla R_m(\Sigma)$.  As a consequence of those results, we deduce that $\nabla R_m(\Sigma)$ is well-defined.

\begin{proposition}
\label{prop_termwise_diff}
Suppose that $\|\Sigma\| \le \gamma_0 d^{r/2}$, where $\gamma_0 > 0$ and $0 \le r < 1$.  Then, for all $m \ge 2$ and all $d \ge (2\gamma_2^2)^{1/(1-r)}$, 
\begin{align}
\label{eq_deriv_Rm_bound}
\|\nabla R_m(\Sigma)\| \le \frac{(2e)^{1/2} \, \gamma_2^{m-1}}{[(m-1)!]^{1/2}} d^{-[1+(m-1)(1-r)]/2}.
\end{align}
In particular, 
$\|\nabla R_m(\Sigma)\| = O\big(d^{-[1 + (m-1)(1-r)]/2}\big)$ as $d \to \infty$.
\end{proposition}

As a consequence of \eqref{eq_deriv_Rm_bound}, the accuracy of the approximation \eqref{eq_matrix_approx} increases as $d \to \infty$ (with $m$ held fixed) or as $m$ increases (with $d$ held fixed).  Also, as with \eqref{eq_Rm_bound_4}, the bound in \eqref{eq_deriv_Rm_bound} increases as $\gamma_0$ or $r$ increases.  

We now obtain an asymptotic expansion of $\Cov(X)$ up to terms of arbitrary degree in $\Sigma$.  

\begin{theorem}
\label{th_cov_matrix}
Suppose that the random vector $X$ has a Bingham distribution with the probability density function \eqref{eq_Bingham_pdf}.  Suppose also that $\|\Sigma\| \le \gamma_0 d^{r/2}$, where $\gamma_0 > 0$ and $0 \le r < 1$.  Then for all $l, m \ge 2$ and all $d > (6\gamma_2^2)^{1/(1-r)}$, 
\begin{equation}
\label{eq_cov_matrix_2}
\Cov(X) = \bigg(1 - \sum_{j=1}^{l-1} \frac{(1/2)_j}{(d/2)_j} \frac{C_{(j)}(\Sigma)}{j!}\bigg)
\bigg(\sum_{k=1}^{m-1} \frac{(1/2)_k}{(d/2)_k} \frac{\nabla C_{(k)}(\Sigma)}{k!}\bigg) + O(d^{-\alpha}),
\end{equation}
where 
$$
\alpha = \tfrac12 \min\{3-2r,1+(m-1)(1-r)\} = 
\begin{cases}
\tfrac12(2-r), & m=2 \\
\tfrac12(3-2r), & m \ge 3
\end{cases}.
$$
\end{theorem}

As we showed before, $C_{(j)}(\Sigma)$ and $\nabla C_{(k)}(\Sigma)$ are homogeneous polynomials of degrees $j-1$ and $k-1$, respectively.  Therefore the highest degree of $\Sigma$ appearing in \eqref{eq_cov_matrix_2} equals $l-1 + m-2 = l+m-3$.  

Consider the asymptotic expansion of $\Cov(X)$ arising from \eqref{eq_cov_matrix_2} with $l = 2$ and $m = 3$.  Applying \eqref{eq_deriv_power_trace}, the rules for derivatives of polynomials in $\Sigma$, to \eqref{eq_zonal_polys}, the formulas for the zonal polynomials $C_{(1)}(\Sigma)$ and $C_{(2)}(\Sigma)$, we have 
$$
\nabla C_{(1)}(\Sigma) = \nabla \tr \Sigma = I_d
$$
and 
$$
\nabla C_{(2)}(\Sigma) = \frac{1}{3} \nabla [(\tr \Sigma)^2 + 2 \tr(\Sigma^2)] = \frac{2}{3} [(\tr \Sigma) I_d + 2 \Sigma].
$$
Substituting these expressions into \eqref{eq_cov_matrix_2} and collecting terms we obtain 
\begin{align*}
\Cov(X) &= \Big(1 - \frac{1}{d} \tr(\Sigma)\Big) \Big(\frac{1}{d} I_d + \frac{1}{d(d+2)} [(\tr \Sigma) I_d + 2 \Sigma]\Big) + O(d^{-(3-2r)/2}) \\
&= \frac{1}{d} I_d - \frac{2}{(d+2)} (\tr \Sigma) I_d + \frac{2}{d(d+2)} \Sigma \\
& \quad - \frac{1}{d^2(d+2)} (\tr\Sigma)^2 I_d - \frac{2}{d^2(d+2)} (\tr\Sigma)\Sigma + O(d^{-(3-2r)/2}).
\end{align*}

In some articles it was assumed, without loss of generality, that $\tr(\Sigma) = 0$; cf.,Kume and Walker \cite{KumeW}.  In that case, the above approximation reduces to 
$$
\Cov(X) = \frac{1}{d} I_d + \frac{2}{d(d+2)} \Sigma + O(d^{-(3-2r)/2}).
$$

\section{Implications for applications}
\label{sec_implications}

The results derived in Sections \ref{sec_asymp_normaliz_const} and \ref{sec_asymp_cov_matrix} provide explicit bounds on the remainder terms in asymptotic expansions, for large values of $d$, of the normalizing constant and the covariance matrix of the Bingham distribution on $\mcalS^{d-1}$.  These results are, to the best of our knowledge, the first such derivations in the literature, so we now comment on the implications of those bounds for real-world applications.  

To assess the implications of our results for applications, we surveyed the literature on the Bingham distributions in a search for the range of reported values of $d$.  Dryden \cite{Dryden} applied the high-dimensional Bingham distribution with $d = 62,501$ to model cortical surfaces using magnetic resonance images of the human brain.  A related study by Brignell, \textit{et al.} \cite{brignell2010surface} applied principal components analysis to reduce the data to dimensions up to $20$.  

Sra \cite{sra2018directional} reviewed applications of high-dimensional spherical data analyses in machine learning.  Those areas include text clustering, gene expression data analysis, feature extraction, and wireless communications, with dimensions up to $d = 1,000$.  Brombin, Pesarin, and Salmaso \cite{brombin2011dealing} discussed examples in the study of shape modeling, and simulations for which the dimension of the data is as large as $d = 50$.  

Dai, Dorman, Dutta, and Maitra \cite{dai2021exploratory} applied spherical data analysis to the characterization of changes in cerebral blood-flow when the human brain is not subject to any stimuli or tasks, to the assessment of variability in handwritten numerical digits, and to the identification of genetic pathways that underlie cancer.  Those topics involve dimensions ranging from $d = 20$ to $d = 100$.  

To summarize, the articles that we reviewed described applications with dimensions ranging from small, $d = 20$, to very large, $d = 62,501$.

\begin{table}[!t]
\caption{Bounds for $\vert R_m(\Sigma)\vert$ and $\|\nabla R_m(\Sigma)\|$ with $(\gamma_0,r) =  (1,0.5)$.}
\label{table_moderate_r}
\begin{subtable}[t]{0.48\textwidth}
\caption{Bounds for $\vert R_m(\Sigma)\vert$}
\qquad
\begin{tabular}[t]{@{} r c *{3}{d{1.5}} @{}}
\toprule
$d$ & \mc{$m=3$} & \mc{$m=6$} & \mc{$m=10$} \\
\midrule
20 &      0.18782 &       0.00349 &       0.00001 \\
25 &      0.15887 &       0.00250 &       0.00000 \\
50 &      0.09447 &       0.00088 &       0.00000 \\
75 &      0.06970 &       0.00048 &       0.00000 \\
100 &     0.05617 &       0.00031 &       0.00000 \\
250 &     0.02825 &       0.00008 &       0.00000 \\
500 &     0.01680 &       0.00003 &       0.00000 \\
750 &     0.01239 &       0.00002 &       0.00000 \\
1000 &    0.00999 &       0.00001 &       0.00000 \\
5000 &    0.00299 &       0.00000 &       0.00000 \\
10000 &   0.00178 &       0.00000 &       0.00000 \\
25000 &   0.00089 &       0.00000 &       0.00000 \\
50000 &   0.00053 &       0.00000 &       0.00000 \\
62501 &   0.00045 &       0.00000 &       0.00000 \\
\bottomrule
\end{tabular}
\label{table_moderate_r_a}
\end{subtable}
\hspace{\fill}
\begin{subtable}[t]{0.48\textwidth}
\caption{Bounds for $\|\nabla R_m(\Sigma)\|$}
\begin{tabular}[t]{@{} r c *{3}{d{1.5}} @{}}
\toprule
$d$ & \mc{$m=3$} & \mc{$m=6$} & \mc{$m=10$} \\
\midrule
20 &      0.15383 &       0.00535 &       0.00002 \\
25 &      0.12306 &       0.00362 &       0.00001 \\
50 &      0.06153 &       0.00108 &       0.00000 \\
75 &      0.04102 &       0.00053 &       0.00000 \\
100 &     0.03077 &       0.00032 &       0.00000 \\
250 &     0.01231 &       0.00006 &       0.00000 \\
500 &     0.00615 &       0.00002 &       0.00000 \\
750 &     0.00410 &       0.00001 &       0.00000 \\
1000 &    0.00308 &       0.00001 &       0.00000 \\
5000 &    0.00062 &       0.00000 &       0.00000 \\
10000 &   0.00031 &       0.00000 &       0.00000 \\
25000 &   0.00012 &       0.00000 &       0.00000 \\
50000 &   0.00006 &       0.00000 &       0.00000 \\
62501 &   0.00005 &       0.00000 &       0.00000 \\
\bottomrule
\end{tabular}
\label{table_moderate_r_b}
\end{subtable}
\end{table}

Consider the bound in \eqref{eq_Rm_bound_4} for $R_m(\Sigma)$, the remainder term in the expansion of $\Psi_d(\Sigma)$.  With $(\gamma_0,r) = (1,0.5)$ and $m = 3, 6, 10$ we computed, for $d > (2\gamma_2^2)^{1/(1-r)} \simeq 14$, that bound.  The computed values, given in Table \ref{table_moderate_r_a}, illustrate that for moderate values of $r$, the bound \eqref{eq_Rm_bound_4} decreases quickly as $m$ increases.  In particular, for $m \ge 6$, all entries in Table \ref{table_moderate_r_a} are less than $10^{-2}$.  

We also computed values of the bound in \eqref{eq_deriv_Rm_bound} for $\|\nabla R_m(\Sigma)\|$.  Several such values are presented in Table \ref{table_moderate_r_b}, and they also illustrate that, for moderate values of $r$, the values of $\|\nabla R_m(\Sigma)\|$ can decrease quickly as $d$ increases.

\begin{table}[!t]
\caption{Bounds for $\vert R_m(\Sigma)\vert$ and $\|\nabla R_m(\Sigma)\|$ with $(\gamma_0,r) = (1,0.75)$.}
\label{table_large_r}
\centering
\begin{subtable}[t]{0.48\textwidth}
\caption{Bounds for $\vert R_m(\Sigma)\vert$}
\qquad
\begin{tabular}[t]{@{} r c *{3}{d{1.5}} @{}}
\toprule
$d$ & \mc{$m=3$} & \mc{$m=6$} & \mc{$m=10$} \\
\midrule
200 &     0.24357 &       0.00587	 &       0.00002	\\
225 &     0.23304 &       0.00538	 &       0.00001	\\
250 &     0.22401 &       0.00497	 &       0.00001	\\
275 &     0.21615 &       0.00463	 &       0.00001	\\
500 &     0.17274 &       0.00295	 &       0.00001	\\
750 &     0.14837 &       0.00218	 &       0.00000	\\
1000	 &    0.13320 &       0.00176	 &       0.00000	\\
2000	 &    0.10271 &       0.00104	 &       0.00000	\\
2500	 &    0.09447 &       0.00088	 &       0.00000	\\
5000	 &    0.07284 &       0.00053	 &       0.00000	\\
10000 &   0.05617 &       0.00031	 &       0.00000	\\
25000 &   0.03984 &       0.00016	 &       0.00000	\\
50000 &   0.03072 &       0.00009	 &       0.00000	\\
62501 &   0.02825 &       0.00008	 &       0.00000	\\
\bottomrule
\end{tabular}
\label{table_large_r_a}
\end{subtable}
\hspace{\fill}
\begin{subtable}[t]{0.48\textwidth}
\caption{Bounds for $\|\nabla R_m(\Sigma)\|$}
\begin{tabular}[t]{@{} r c *{3}{d{1.5}} @{}}
\toprule
$d$ & \mc{$m=3$} & \mc{$m=6$} & \mc{$m=10$} \\
\midrule
200 &       0.05785 &     0.00261 &     0.00001 \\
225 &       0.05296 &     0.00229 &     0.00001 \\
250 &       0.04893 &     0.00203 &     0.00001 \\
275 &       0.04556 &     0.00182 &     0.00001 \\
500 &       0.02910 &     0.00093 &     0.00000 \\
750 &       0.02147 &     0.00059 &     0.00000 \\
1000 &      0.01730 &     0.00043 &     0.00000 \\
2000 &      0.01029 &     0.00020 &     0.00000 \\
2500 &      0.00870 &     0.00015 &     0.00000 \\
5000 &      0.00517 &     0.00007 &     0.00000 \\
10000 &     0.00308 &     0.00003 &     0.00000 \\
25000 &     0.00155 &     0.00001 &     0.00000 \\
50000 &     0.00092 &     0.00001 &     0.00000 \\
62501 &     0.00078 &     0.00000 &     0.00000 \\
\bottomrule
\end{tabular}
\label{table_large_r_b}
\end{subtable}
\end{table}

As we noted before, the upper bounds in \eqref{eq_Rm_bound_4} and \eqref{eq_deriv_Rm_bound} are increasing functions of $r$.  Therefore the upper bounds derived from any $r < 0.5$ will be less than the corresponding values in Table \ref{table_moderate_r}.  

As regards values of $r$ greater than $0.5$ we note that if $r$ is very large, e.g., $r > 0.9$, the restriction $d > (2\gamma_2^2)^{1/(1-r)}$ can result in extremely large values of $d$, even in the tens of millions, which have never arisen in real-world problems.  Thus, we calculated the upper bounds in \eqref{eq_Rm_bound_4} and \eqref{eq_deriv_Rm_bound} for $(\gamma_0,r) = (1,0.75)$, $m = 3, 6, 10$, and $d > (2\gamma_2^2)^{1/(1-r)} \simeq 194$.  These computed values are provided in Table \ref{table_large_r}, and they also illustrate the rapid decrease in the values of $\vert R_m(\Sigma)\vert$ and $\|\nabla R_m(\Sigma)\|$ as $d$ increases.  

It is noticeable that, in Table \ref{table_moderate_r_a} with $m=3$, each entry is greater than the corresponding entry in Table \ref{table_moderate_r_b}, and the order is reversed for $m = 6$ and $m = 10$.  Throughout Table \ref{table_large_r_a}, however, every entry is greater than or equal to the corresponding entry in Table \ref{table_large_r_b}.   This phenomenon may be explained by comparing the upper bounds on the right-hand sides of \eqref{eq_Rm_bound_4} and \eqref{eq_deriv_Rm_bound}; on doing so, we determine that a necessary and sufficient condition for the bound in \eqref{eq_Rm_bound_4} to be less than or equal to the bound in \eqref{eq_deriv_Rm_bound} is that 
\begin{equation}
\label{eq_compare_tables}
d^r \ge \frac{m(m+1)}{4\gamma_0}.
\end{equation}
Consider Table \ref{table_moderate_r}, where we have $(\gamma_0,r) = (1,0.5)$.  Then the condition \eqref{eq_compare_tables} is satisfied for $m = 3$ and all $d \ge 9$; this explains why each entry in the first column of Table \ref{table_moderate_r_a} is greater than the corresponding entry in Table \ref{table_moderate_r_b}.  For $m = 6$ the condition \eqref{eq_compare_tables} is satisfied in Table \ref{table_moderate_r} for $d \ge 111$, and this is illustrated by comparing the entries of the second column of that table.  Similar remarks can be made about Table \ref{table_large_r} and, more generally, the condition \eqref{eq_compare_tables} implies that, for all sufficiently large $d$, the bound provided by \eqref{eq_deriv_Rm_bound} will be less than the corresponding bound from \eqref{eq_Rm_bound_4}.  

We conclude this section by suggesting a method for choosing $m$.  A statistician undertaking an application of the high-dimensional Bingham distribution first should conduct exploratory work to determine plausible values for $\gamma_0$ and $r$.  Next, the statistician should choose a tolerance or cutoff level, denoted by $\epsilon$, for the larger of the upper bounds in Propositions \ref{prop_Ck_series_tail_bound} and \ref{prop_termwise_diff}; the value of $\epsilon$ should reflect the importance of the application, so that $\epsilon$ presumably would be chosen to be small in the case of crucial biomedical research.  Once $\gamma_0$, $r$, and $\epsilon$ have been chosen then $m$ should be chosen so that the larger of those two upper bounds is no greater than $\epsilon$, i.e., 
$$
\max\left\{\frac{\gamma_3 \gamma_2^m}{[(m+1)!]^{1/2}} d^{-m(1-r)/2},\frac{(2e)^{1/2} \, \gamma_2^{m-1}}{[(m-1)!]^{1/2}} d^{-[1+(m-1)(1-r)]/2}\right\} \le \epsilon.
$$

\bigskip
\bigskip

\noindent
\textbf{Acknowledgments}.  The authors are grateful to the referees and the editors for helpful comments on the initial version of this article.

\medskip

\noindent
\textbf{Declarations}.  No funds, grants, or other support were received for conducting this research.  The authors have no relevant financial or non-financial interests to disclose.  The authors have no conflicts of interest to declare.

\appendix

\section{Appendix: Proofs}
\label{sec_appendix}

\subsection{The proofs of Proposition \ref{prop_Ck_series_tail_bound} and Corollary \ref{cor_inverse_normaliz_const}}

\noindent
\textit{Proof of Proposition \ref{prop_Ck_series_tail_bound}}: 
By \eqref{eq_remainder_term}, 
\begin{equation}
\label{eq_Rm_bound}
\vert R_m(\Sigma) \vert \le \sum_{k=m}^\infty \frac{(1/2)_k}{(d/2)_k} \frac{\vert C_{(k)}(\Sigma) \vert}{k!},
\end{equation}
and by \eqref{eq_Ck_multisum}, 
\begin{equation}
\label{eq_Ck_multisum_abs}
\vert C_{(k)}(\Sigma) \vert \le \frac{k!}{(1/2)_k} \ \sum_{i_1 + 2 i_2 + 3 i_3 + \cdots + k i_k = k} \ \prod_{j=1}^k \frac{\vert \tr (\Sigma^j) \vert^{i_j}}{i_j! \, (2j)^{i_j}}.
\end{equation}
For each vector of indices $(i_1,i_2,\ldots,i_k)$ such that $i_1 + 2 i_2 + \cdots + k i_k = k$, we have  
\begin{align}
\label{eq_reznick_1}
\prod_{j=1}^k \vert \tr (\Sigma^j) \vert^{i_j} &= \|\Sigma\|^k \prod_{j=1}^k \frac{\vert\tr (\Sigma^j)\vert^{i_j}}{\|\Sigma\|^{ji_j}} \nonumber \\
&= \|\Sigma\|^k \prod_{j=1}^k \frac{\vert\tr (\Sigma^j)\vert^{i_j}}{(\tr (\Sigma^2))^{ji_j/2}}.
\end{align}
Denote by $\lambda_1,\ldots,\lambda_d$ the eigenvalues of $\Sigma$.  Then for each $j=1,\ldots,k$, 
$$
\vert\tr (\Sigma^j)\vert^{i_j} = \bigg\vert\sum_{i=1}^d \lambda_i^j\bigg\vert^{i_j} \le \bigg(\sum_{i=1}^d \vert\lambda_i\vert^j\bigg)^{i_j},
$$
and by substituting this bound into \eqref{eq_reznick_1} we obtain 
\begin{equation}
\label{eq_trace_eigen}
\prod_{j=1}^k \vert \tr (\Sigma^j) \vert^{i_j} \le \|\Sigma\|^k \prod_{j=1}^k \frac{\left(\sum_{i=1}^d \vert\lambda_i\vert^j\right)^{i_j}}{\left(\sum_{i=1}^d \vert\lambda_i\vert^2\right)^{ji_j/2}}.
\end{equation}
By a remarkable result of Reznick \cite[p.~447, eq.~(2.12)]{Reznick} we obtain, for each $j=1,\ldots,k$ and all $d$,  
\begin{equation}
\label{eq_reznick_3}
\frac{\left(\sum_{i=1}^d \vert\lambda_i\vert^j\right)^{i_j}}{\left(\sum_{i=1}^d \vert\lambda_i\vert^2\right)^{ji_j/2}} \le d^{\max\{0,(2-j) i_j/2\}} = 
\begin{cases}
d^{i_1/2}, & j=1 \\
1, & j \ge 2
\end{cases}.
\end{equation}
Substituting the bound in \eqref{eq_reznick_3} into \eqref{eq_trace_eigen}, we obtain 
\begin{equation}
\label{eq_reznick_4}
\prod_{j=1}^k \vert\tr (\Sigma^j)\vert^{i_j} \le d^{i_1/2} \|\Sigma\|^k \equiv d^{i_1/2} \|\Sigma\|^{i_1+2i_2+3i_3+\cdots+ki_k}.
\end{equation}
On applying \eqref{eq_reznick_4} to \eqref{eq_Ck_multisum_abs} we obtain 
\begin{equation}
\label{eq_Ck_bound_2}
\vert C_{(k)}(\Sigma)\vert \le \frac{k!}{(1/2)_k} a_k \|\Sigma\|^k,
\end{equation}
where 
\begin{equation}
\label{eq_ak}
a_k = \sum_{i_1 + 2 i_2 + 3 i_3 + \cdots + k i_k = k} d^{i_1/2} \prod_{j=1}^k \frac{1}{i_j! \, (2j)^{i_j}}.
\end{equation}
We will derive in Lemma \ref{lem_ak_formula} a single-sum expression for the multiple sum $a_k$, and then we will deduce from that single-sum the inequality 
\begin{equation}
\label{eq_ak_bound}
a_k \le \frac{(d^{1/2}/2)_k}{k!}.
\end{equation}
Substituting \eqref{eq_ak_bound} into \eqref{eq_Ck_bound_2} and recalling that $\|\Sigma\| \le \gamma_0 d^{r/2}$ for all $d$, we obtain 
\begin{equation}
\label{eq_Ck_bound_large_d}
\vert C_{(k)}(\Sigma)\vert \le \frac{(d^{1/2}/2)_k}{(1/2)_k} \|\Sigma\|^k \le \frac{(d^{1/2}/2)_k}{(1/2)_k} (\gamma_0 d^{r/2})^k.
\end{equation}
Inserting \eqref{eq_Ck_bound_large_d} into \eqref{eq_Rm_bound} we obtain 
\begin{equation}
\label{eq_Rm_bound_2}
\vert R_m(\Sigma) \vert \le \sum_{k=m}^\infty \frac{(d^{1/2}/2)_k}{(d/2)_k} \frac{\gamma_0 ^k d^{rk/2}}{k!}.
\end{equation}

As we shall show in Lemma \ref{lem_rem_rk_ineq}, with $\gamma_1 = \tfrac12(\sqrt{3}+1) \simeq 1.366025$, 
\begin{equation}
\label{eq_rk_bound_3}
\frac{(d^{1/2}/2)_k}{(d/2)_k} \le \gamma_1^{k-1} ((k-1)!)^{1/2} d^{-k/2}
\end{equation}
for all $d \ge 1$, $k \ge 1$.  Letting $\gamma_2 = \gamma_1 \gamma_0$ and substituting \eqref{eq_rk_bound_3} into \eqref{eq_Rm_bound_2}, we find that for all $d$, 
\begin{equation}
\label{eq_Rm_bound_31}
\vert R_m(\Sigma)\vert \le \gamma_1^{-1} \sum_{k=m}^\infty \frac{((k-1)!)^{1/2}}{k!} \, \gamma_2^k d^{-k(1-r)/2}.
\end{equation}
By applying the ratio test, we find that the series \eqref{eq_Rm_bound_31} converges absolutely for all fixed $d$.  

Note that although the $k$th term in \eqref{eq_Rm_bound_31} is $O(d^{-k(1-r)/2})$, it is not evident that the sum of the resulting series is $O(d^{-m(1-r)/2})$ as $d \to \infty$.  To establish that stated convergence rate we apply to \eqref{eq_Rm_bound_31} the Cauchy-Schwarz inequality, obtaining 
\begin{align}
\label{eq_Rm_bound_3}
\vert R_m(\Sigma) \vert &\le \gamma_1^{-1} \Big(\sum_{k=m}^\infty \frac{(k-1)!}{(k!)^2}\Big)^{1/2} \Big(\sum_{k=m}^\infty \gamma_2^{2k} d^{-k(1-r)}\Big)^{1/2}
\nonumber \\
&= \gamma_1^{-1} \Big(\sum_{k=m}^\infty \frac{1}{k \cdot k!}\Big)^{1/2} \cdot \gamma_2^m d^{-m(1-r)/2}  \big(1 - \gamma_2^2 d^{-(1-r)}\big)^{-1/2}.
\end{align}
For $k \ge m$, it is straightforward that 
$$
\frac{m}{k} \le 1 \le \binom{k}{m} \equiv \frac{k!}{m! (k-m)!},
$$
equivalently, 
$$
\frac{1}{k \cdot k!} \le \frac{1}{m \cdot m!} \cdot \frac{1}{(k-m)!},
$$
also that 
$$
\frac{1}{m} \le \frac{2}{m+1}.
$$
Therefore 
$$
\frac{1}{k \cdot k!} \le \frac{1}{m \cdot m!} \cdot \frac{1}{(k-m)!} \le \frac{2}{(m+1)!} \cdot \frac{1}{(k-m)!}.
$$
Summing this inequality over all $k \ge m$ we obtain 
\begin{equation}
\label{eq_Ei_bound}
\sum_{k=m}^\infty \frac{1}{k \cdot k!} \le \frac{2}{(m+1)!} \sum_{k=m}^\infty \frac{1}{(k-m)!} = \frac{2e}{(m+1)!}.
\end{equation}
Noting that the inequality $d \ge (2\gamma_2^2)^{1/(1-r)}$ is equivalent to 
\begin{equation}
\label{eq_Rm_CS_bound}
\big(1 - \gamma_2^2 d^{-(1-r)}\big)^{-1} \le 2,
\end{equation}
and applying \eqref{eq_Ei_bound} and \eqref{eq_Rm_CS_bound} to \eqref{eq_Rm_bound_3} we obtain, for all $d \ge (2\gamma_2^2)^{1/(1-r)}$, the inequality  
\begin{equation}
\label{eq_Rm_bound_45}
\vert R_m(\Sigma) \vert \le 2 \gamma_1^{-1} \Big(\frac{2e}{(m+1)!}\Big)^{1/2} \gamma_2^m d^{-m(1-r)/2}.
\end{equation}
This establishes \eqref{eq_Rm_bound_4} and proves that $R_m(\Sigma) = O(d^{-m(1-r)/2})$ as $d \to \infty$, and the proof of the proposition now is complete.  
$\Box$

\bigskip

\noindent
\textit{Proof of Corollary \ref{cor_inverse_normaliz_const}}: 
By \eqref{eq_Rm_bound_45} with $m=1$ we have, for all $d \ge (2\gamma_2^2)^{1/(1-r)}$, 
\begin{equation}
\label{eq_R1_bound}
\vert R_1(\Sigma) \vert \le 2 \gamma_1^{-1} e^{1/2} \gamma_2 d^{-(1-r)/2}.
\end{equation}
For $d > (6\gamma_2^2)^{1/(1-r)}$, we use the values of the constants given in \eqref{eq_constants} to determine that the right-hand side of \eqref{eq_R1_bound} is strictly less than $1$.  Then we apply the geometric series to obtain, for all $l \ge 2$, 
\begin{align}
\label{eq_Psi_inv_exp}
[\Psi_d(\Sigma)]^{-1} &= [1 + R_1(\Sigma)]^{-1} \nonumber \\
&= 1 - R_1(\Sigma) + \sum_{j=2}^\infty (-1)^j \big(R_1(\Sigma)\big)^j \nonumber \\
&= 1 - \bigg(\sum_{j=1}^{l-1} \frac{(1/2)_j}{(d/2)_j} \frac{C_{(j)}(\Sigma)}{j!} + R_l(\Sigma)\bigg) + \sum_{j=2}^\infty (-1)^j \big(R_1(\Sigma)\big)^j.
\end{align}
By Proposition \ref{prop_Ck_series_tail_bound}, $R_l(\Sigma) = O(d^{-l(1-r)/2})$.  On applying \eqref{eq_Rm_bound_45} with $m=1$ we obtain, for all $d > (6\gamma_2^2)^{1/(1-r)}$, 
\begin{align}
\label{eq_Psi_inv_bound}
\left\vert\sum_{j=2}^\infty (-1)^j \big(R_1(\Sigma)\big)^j\right\vert &\le \sum_{j=2}^\infty \vert R_1(\Sigma)\vert^j \nonumber \\
&\le \sum_{j=2}^\infty \big[2 \gamma_1^{-1} \gamma_2 e^{1/2} d^{-(1-r)/2}\big]^j.
\end{align}
By summing this geometric series, we deduce that \eqref{eq_Psi_inv_bound} is $O(d^{-(1-r)})$ as $d \to \infty$.  Therefore by \eqref{eq_Psi_inv_exp}, for all $l \ge 2$ and all $d > (6\gamma_2^2)^{1/(1-r)}$, 
$$
[\Psi_d(\Sigma)]^{-1} = 1 - \sum_{j=1}^{l-1} \frac{(1/2)_j}{(d/2)_j} \frac{C_{(j)}(\Sigma)}{j!} + O(d^{-l(1-r)/2}) + O(d^{-(1-r)}).
$$
On applying the well-known property, 
\begin{equation}
\label{eq_sum_of_Big_Os}
O(d^{-a}) + O(d^{-b}) = O(d^{-\min(a,b)})
\end{equation}
for $a,b > 0$, we obtain \eqref{eq_Psi_d_inverse}.  
$\Box$

\subsection{The proof of Proposition \ref{prop_termwise_diff} and Theorem \ref{th_cov_matrix}}

\noindent
{\it Proof of Lemma \ref{lem_nabla_app}}. \ 
The result \eqref{eq_deriv_exp_trace} follows by straightforwardly applying each element of the matrix $\nabla$ to the function $\exp\big(\tr (\Sigma H)\big)$.  

The formula \eqref{eq_deriv_power_trace} follows from the chain rule: 
$$
\nabla [\tr (\Sigma)]^k = k [\tr (\Sigma)]^{k-1} \nabla \tr (\Sigma) = k [\tr (\Sigma)]^{k-1} I_d.
$$

As for \eqref{eq_deriv_trace_power}, that result can be deduced by arguments similar to those given by Sebastiani \cite[Lemmas 3.1 or 5.1]{Sebastiani}.  Also see Dwyer and Macphail \cite[p.~528, Section 14]{Dwyer} for the derivatives of the trace of powers of square, \textit{non-symmetric} matrices; calculations similar to theirs can also lead to \eqref{eq_deriv_trace_power}.  

A succinct derivation of \eqref{eq_deriv_trace_power} is obtained from the Taylor expansion \eqref{eq_Taylor_expansion}, as follows.  By using the commutativity property of the trace, i.e., $\tr(\Sigma H) = \tr(H \Sigma)$, we find that 
\begin{align}
\label{eq_exp_trace_power}
\tr [(\Sigma+H)^k] &= \tr[(\Sigma+H)(\Sigma+H)\cdots(\Sigma+H)] \nonumber \\
&= \tr(\Sigma^k) + k \tr(H \Sigma^{k-1}) + O(\|H\|^2) \nonumber \\
&\equiv \tr(\Sigma^k) + \langle H,k\Sigma^{k-1}\rangle + O(\|H\|^2)
\end{align}
as $H \to 0$.  Setting $f(\Sigma) = \tr(\Sigma^k)$ in \eqref{eq_Taylor_expansion} and comparing the result with \eqref{eq_exp_trace_power}, we obtain \eqref{eq_deriv_trace_power}.  
$\Box$

\noindent
\textit{Proof of Proposition \ref{prop_termwise_diff}}: 
By the subadditivity property of the Frobenius norm, we have 
\begin{align*}
\|\nabla R_m(\Sigma)\| &= \bigg\|\sum_{k=m}^\infty \frac{(1/2)_k}{(d/2)_k} \frac{\nabla C_{(k)}(\Sigma)}{k!}\bigg\| \\
&\le \sum_{k=m}^\infty \bigg\|\frac{(1/2)_k}{(d/2)_k} \frac{\nabla C_{(k)}(\Sigma)}{k!}\bigg\|.
\end{align*}
For $k \ge 2$, it follows from \eqref{eq_Ck_multisum} that 
\begin{align}
\label{eq_Ck_pk}
C_{(k)}(\Sigma) &= \frac{k!}{(1/2)_k} \Bigg[\frac{[\tr (\Sigma)]^k}{k! \, 2^k} + \sum_{\substack{i_1 + 2 i_2 + 3 i_3 + \cdots + k i_k = k \\ i_1 \le k-2}} \ \prod_{j=1}^k \frac{[\tr (\Sigma^j)]^{i_j}}{i_j! \, (2j)^{i_j}}\Bigg] \nonumber \\
&\equiv \frac{[\tr (\Sigma)]^k}{(1/2)_k \, 2^k} + \frac{k!}{(1/2)_k} p_k(\Sigma),
\end{align}
where 
\begin{equation}
\label{eq_pk}
p_k(\Sigma) = \sum_{\substack{i_1 + 2 i_2 + 3 i_3 + \cdots + k i_k = k \\ i_1 \le k-2}} \ \prod_{j=1}^k \frac{[\tr (\Sigma^j)]^{i_j}}{i_j! \, (2j)^{i_j}}
\end{equation}
is a polynomial in $\{\tr(\Sigma),\tr(\Sigma^2),\ldots,\tr(\Sigma^k)\}$.  It follows from \eqref{eq_pk} that $p_k(\Sigma)$ is homogeneous of degree $k$ in $\Sigma$ and its coefficients do not depend on $d$; moreover, because of the restriction $i_1 \le k-2$,d the highest power of $\tr(\Sigma)$ which can appear in $p_k(\Sigma)$ is $k-2$.

By \eqref{eq_deriv_power_trace} and \eqref{eq_Ck_pk}, 
\begin{align}
\label{eq_Ck_derivative}
\frac{(1/2)_k}{(d/2)_k} \frac{\nabla C_{(k)}(\Sigma)}{k!} &= \frac{(1/2)_k}{(d/2)_k \, k!} \nabla \bigg[\frac{[\tr (\Sigma)]^k}{(1/2)_k \, 2^k} + \frac{k!}{(1/2)_k} p_k(\Sigma)\bigg] \nonumber \\
&= \frac{1}{(d/2)_k}\bigg[\frac{1}{(k-1)! \, 2^k} [\tr(\Sigma)]^{k-1} I_d + \nabla p_k(\Sigma)\bigg].
\end{align}
On applying to \eqref{eq_Ck_derivative} the subadditivity property of the Frobenius norm we obtain 
\begin{align}
\label{eq_Ck_norm}
\bigg\|\frac{(1/2)_k}{(d/2)_k} & \frac{\nabla C_{(k)}(\Sigma)}{k!}\bigg\| \nonumber \\
&\le \frac{1}{(d/2)_k} \bigg[\frac{1}{(k-1)! \, 2^k} \vert\tr (\Sigma)\vert^{k-1} \|I_d\| + \|\nabla p_k(\Sigma)\|\bigg] \nonumber \\
&= \frac{1}{(d/2)_k} \bigg[\frac{d^{1/2}}{(k-1)! \, 2^k} \vert\tr (\Sigma)\vert^{k-1} + \|\nabla p_k(\Sigma)\|\bigg].
\end{align}

By \eqref{eq_pk} and the product rule for derivatives, 
\begin{align*}
\nabla p_k(\Sigma) 
&= \sum_{\substack{i_1 + 2 i_2 + 3 i_3 + \cdots + k i_k = k \\ i_1 \le k-2}} \, \nabla \prod_{j=1}^k \frac{[\tr (\Sigma^j)]^{i_j}}{i_j! \, (2j)^{i_j}} \\
&= \sum_{\substack{i_1 + 2 i_2 + 3 i_3 + \cdots + k i_k = k \\ i_1 \le k-2}} \, \sum_{l=1}^k \Big(\prod_{\substack{j=1 \\ j \neq l}}^k \frac{[\tr (\Sigma^j)]^{i_j}}{i_j! \, (2j)^{i_j}}\Big) \cdot \frac{\nabla [\tr (\Sigma^l)]^{i_l}}{i_l! \, (2l)^{i_l}}.
\end{align*}
By \eqref{eq_deriv_trace_power} and the chain rule, 
$$
\nabla [\tr (\Sigma^l)]^{i_l} = i_l [\tr (\Sigma^l)]^{i_l - 1} \nabla [\tr (\Sigma^l)] = l \, i_l [\tr (\Sigma^l)]^{i_l - 1} \Sigma^{l-1};
$$
hence, 
$$
\nabla p_k(\Sigma) 
= \sum_{\substack{i_1 + 2 i_2 + \cdots + k i_k = k \\ i_1 \le k-2}} \sum_{l=1}^k \frac{l \, i_l}{i_l! (2l)^{i_l}} \Big(\prod_{\substack{j=1 \\ j \neq l}}^k \frac{[\tr (\Sigma^j)]^{i_j}}{i_j! (2j)^{i_j}}\Big) [\tr (\Sigma^l)]^{i_l - 1} \Sigma^{l-1},
$$
which also reveals that $\nabla p_k(\Sigma)$ is homogeneous of degree $k-1$.  On applying the subadditivity property of the Frobenius norm, we obtain 
\begin{align*}
& \| \nabla p_k(\Sigma)\| \\
& \ \le \sum_{\substack{i_1 + 2 i_2 + \cdots + k i_k = k \\ i_1 \le k-2}} \bigg(\prod_{j=1}^k \frac{1}{i_j! (2j)^{i_j}}\bigg) \sum_{l=1}^k l \, i_l \Big(\prod_{\substack{j=1 \\ j \neq l}}^k \vert\tr (\Sigma^j)\vert^{i_j}\Big) \vert\tr (\Sigma^l)\vert^{i_l - 1} \|\Sigma^{l-1}\|.
\end{align*}
On applying \eqref{eq_reznick_4} we have, for all $d$, 
$$
\Big(\prod_{\substack{j=1 \\ j \neq l}}^k \vert\tr (\Sigma^j)\vert^{i_j}\Big) \vert\tr (\Sigma^l)\vert^{i_l - 1} \le d^{i_1/2} \|\Sigma\|^{i_1+2i_2+3i_3+\cdots+ki_k-l} = d^{i_1/2} \|\Sigma\|^{k-l}.
$$
By also applying the submultiplicative inequality, $\|\Sigma^{l-1}\| \le \|\Sigma\|^{l-1}$, we obtain 
\begin{align}
\label{eq_pk_derivative}
\|\nabla p_k(\Sigma)\| &\le \sum_{\substack{i_1 + 2 i_2 + \cdots + k i_k = k \\ i_1 \le k-2}} \bigg(\prod_{j=1}^k \frac{1}{i_j! \, (2j)^{i_j}}\bigg) \sum_{l=1}^k l i_l \, d^{i_1/2} \|\Sigma\|^{k-l} \|\Sigma\|^{l-1} \nonumber \\
&= \|\Sigma\|^{k-1} \sum_{\substack{i_1 + 2 i_2 + \cdots + k i_k = k \\ i_1 \le k-2}} d^{i_1/2}\bigg(\prod_{j=1}^k \frac{1}{i_j! \, (2j)^{i_j}}\bigg) \sum_{l=1}^k l \, i_l \nonumber \\
&= k \|\Sigma\|^{k-1} \sum_{\substack{i_1 + 2 i_2 + \cdots + k i_k = k \\ i_1 \le k-2}} d^{i_1/2} \prod_{j=1}^k \frac{1}{i_j! \, (2j)^{i_j}}.
\end{align}
By \eqref{eq_ak} and \eqref{eq_ak_bound}, 
\begin{align*}
\sum_{\substack{i_1 + 2 i_2 + \cdots + k i_k = k \\ i_1 \le k-2}} d^{i_1/2} \prod_{j=1}^k \frac{1}{i_j! \, (2j)^{i_j}} &\equiv a_k - \frac{d^{k/2}}{k! \, 2^k} \\
&\le \frac{1}{k!}\big[(d^{1/2}/2)_k - (d^{1/2}/2)^k\big],
\end{align*}
and by substituting this bound into \eqref{eq_pk_derivative}, we obtain 
$$
\|\nabla p_k(\Sigma)\| \le \frac{1}{(k-1)!}\big[(d^{1/2}/2)_k - (d^{1/2}/2)^k\big] \|\Sigma\|^{k-1}.
$$
On applying the latter inequality at \eqref{eq_Ck_norm}, we find that 
\begin{multline*}
\bigg\| \frac{(1/2)_k}{(d/2)_k} \frac{\nabla C_{(k)}(\Sigma)}{k!}\bigg\| \\
\le \frac{1}{(d/2)_k \, (k-1)!} \bigg[\frac{d^{1/2}}{2^k} \vert\tr (\Sigma)\vert^{k-1} + \big[(d^{1/2}/2)_k - (d^{1/2}/2)^k\big] \|\Sigma\|^{k-1}\bigg].
\end{multline*}
Setting $(i_1,i_2,\ldots,i_k) = (1,0,\ldots,0)$ in \eqref{eq_reznick_4}, we obtain 
$$
\vert\tr (\Sigma)\vert \le d^{1/2} [\tr (\Sigma^2)]^{1/2} = d^{1/2} \|\Sigma\|;
$$
hence 
\begin{align*}
\bigg\| & \frac{(1/2)_k}{(d/2)_k} \frac{\nabla C_{(k)}(\Sigma)}{k!}\bigg\| \\
&\le \frac{1}{(d/2)_k \, (k-1)!} \bigg[\frac{d^{1/2}}{2^k} (d^{1/2}\|\Sigma\|)^{k-1} + \big[(d^{1/2}/2)_k - (d^{1/2}/2)^k\big] \|\Sigma\|^{k-1}\bigg] \\
&= \frac{(d^{1/2}/2)_k}{(d/2)_k} \frac{\|\Sigma\|^{k-1}}{(k-1)!}.
\end{align*}
On applying the inequalities $\|\Sigma\| \le \gamma_0 d^{r/2}$ and \eqref{eq_rk_bound_1} we obtain 
\begin{align*}
\bigg\|\frac{(1/2)_k}{(d/2)_k} \frac{\nabla C_{(k)}(\Sigma)}{k!}\bigg\| &\le \frac{(d^{1/2}/2)_k}{(d/2)_k} \frac{(\gamma_0 d^{r/2})^{k-1}}{(k-1)!} \\
&\le \gamma_1^{k-1} ((k-1)!)^{1/2} d^{-k/2} \frac{(\gamma_0 d^{r/2})^{k-1}}{(k-1)!} \\
&= d^{-1/2} \frac{(\gamma_2 d^{-(1-r)/2})^{k-1}}{((k-1)!)^{1/2}}.
\end{align*}
Therefore 
\begin{align}
\label{eq_nabla_Rm_bound}
\|\nabla R_m(\Sigma)\| &\le \sum_{k=m}^\infty \bigg\|\frac{(1/2)_k}{(d/2)_k} \frac{\nabla C_{(k)}(\Sigma)}{k!}\bigg\| \nonumber \\
&\le d^{-1/2} \sum_{k=m}^\infty \frac{(\gamma_2 d^{-(1-r)/2})^{k-1}}{((k-1)!)^{1/2}},
\end{align}
and, by applying the ratio test, we find that the latter series converges for all fixed $d$.  

On applying the Cauchy-Schwarz inequality to \eqref{eq_nabla_Rm_bound}, we obtain 
$$
\|\nabla R_m(\Sigma)\| \le d^{-1/2} \Big(\sum_{k=m}^\infty \frac{1}{(k-1)!}\Big)^{1/2} \Big(\sum_{k=m}^\infty (\gamma_2^2 d^{-(1-r)})^{k-1}\Big)^{1/2}.
$$
For any nonnegative integers $a$ and $b$, it is elementary that 
$$
\frac{1}{(a+b)!} \le \frac{1}{a! b!};
$$
therefore 
$$
\sum_{k=m}^\infty \frac{1}{(k-1)!} = \sum_{k=0}^\infty \frac{1}{(k+m-1)!} \le \sum_{k=0}^\infty \frac{1}{k! \, (m-1)!} = \frac{e}{(m-1)!}.
$$
Also, for all $d$ such that $\gamma_2^2 d^{-(1-r)} \le 1/2$, equivalently, $d \ge (2\gamma_2^2)^{1/(1-r)}$, we have 
\begin{align*}
\sum_{k=m}^\infty (\gamma_2^2 d^{-(1-r)})^{k-1} &= (\gamma_2^2 d^{-(1-r)})^{m-1} (1 - \gamma_2^2 d^{-(1-r)})^{-1} \\
&\le 2 (\gamma_2^2 d^{-(1-r)})^{m-1},
\end{align*}
and then we obtain 
\begin{align*}
\|\nabla R_m(\Sigma)\| &\le d^{-1/2} \Big(\frac{e}{(m-1)!}\Big)^{1/2} \big(2 (\gamma_2^2 d^{-(1-r)})^{m-1}\big)^{1/2} \\
&= (2e)^{1/2} [(m-1)!]^{-1/2} \gamma_2^{m-1} d^{-[1+(m-1)(1-r)]/2}.
\end{align*}
Therefore $\|\nabla R_m(\Sigma)\| = O(d^{-[1+(m-1)(1-r)]/2})$ as $d \to \infty$. 
$\Box$

\bigskip

\noindent
\textit{Proof of Theorem \ref{th_cov_matrix}}: 
Since $\nabla C_{(0)}(\Sigma) = 0$ then by the zonal polynomial expansion \eqref{eq_c_series} of $\Psi_d(\Sigma)$, 
\begin{equation}
\label{eq_nabla_Sigma}
\nabla \Psi_d(\Sigma) = \sum_{k=1}^{m-1} \frac{(1/2)_k}{(d/2)_k} \frac{\nabla C_{(k)}(\Sigma)}{k!} + \nabla R_m(\Sigma)
\end{equation}
and, by Proposition \ref{prop_termwise_diff}, $\|\nabla R_m(\Sigma)\| = O(d^{-[1+(m-1)(1-r)]/2})$ as $d \to \infty$.  

On applying the asymptotic expansion of $[\Psi_d(\Sigma)]^{-1}$ given in \eqref{eq_Psi_d_inverse}, and using \eqref{eq_nabla_Sigma}, we obtain 
\begin{align*}
\Cov(X) &= [\Psi_d(\Sigma)]^{-1} \nabla \Psi_d(\Sigma) \\
&= \bigg(1 - \sum_{j=1}^{l-1} \frac{(1/2)_j}{(d/2)_j} \frac{C_{(j)}(\Sigma)}{j!} + O(d^{-(1-r)})\bigg) \\
&\qquad \cdot \bigg(\sum_{k=1}^{m-1} \frac{(1/2)_k}{(d/2)_k} \frac{\nabla C_{(k)}(\Sigma)}{k!} + O(d^{-[1+(m-1)(1-r)]/2})\bigg).
\end{align*}
Expanding this product, we obtain 
\begin{align}
\label{eq_CovX_BigO}
\Cov(X) &= \bigg(1 - \sum_{j=1}^{l-1} \frac{(1/2)_j}{(d/2)_j} \frac{C_{(j)}(\Sigma)}{j!}\bigg) \bigg(\sum_{k=1}^{m-1} \frac{(1/2)_k}{(d/2)_k} \frac{\nabla C_{(k)}(\Sigma)}{k!}\bigg) \nonumber \\
&\qquad + O(d^{-(1-r)}) \cdot \bigg(\sum_{k=1}^{m-1} \frac{(1/2)_k}{(d/2)_k} \frac{\nabla C_{(k)}(\Sigma)}{k!}\bigg) \nonumber \\
&\qquad + 1 \cdot O(d^{-[1+(m-1)(1-r)]/2}) \nonumber \\
&\qquad + \bigg(\sum_{j=1}^{l-1} \frac{(1/2)_j}{(d/2)_j} \frac{C_{(j)}(\Sigma)}{j!}\bigg) \cdot O(d^{-[1+(m-1)(1-r)]/2}) \nonumber \\
&\qquad + O(d^{-(1-r)}) \cdot O(d^{-[1+(m-1)(1-r)]/2}).
\end{align}
On applying Proposition \ref{prop_termwise_diff}, we obtain 
\begin{align*}
\bigg\|\sum_{k=1}^{m-1} \frac{(1/2)_k}{(d/2)_k} \frac{\nabla C_{(k)}(\Sigma)}{k!}\bigg\| &= \|\nabla R_1(\Sigma) - \nabla R_m(\Sigma)\| \\
&\le \|\nabla R_1(\Sigma)\| + \|\nabla R_m(\Sigma)\| \\
&= O(d^{-1/2}) + O(d^{-[1+(m-1)(1-r)]/2}) \\
&= O(d^{-1/2});
\end{align*}
therefore 
\begin{align*}
O(d^{-(1-r)}) \cdot \bigg\|\sum_{k=1}^{m-1} \frac{(1/2)_k}{(d/2)_k} \frac{\nabla C_{(k)}(\Sigma)}{k!}\bigg\| &= O(d^{-(1-r)}) \cdot O(d^{-1/2}) \\
&= O(d^{-(3-2r)/2}).
\end{align*}
Similarly, by Proposition \ref{prop_Ck_series_tail_bound}, 
\begin{align*}
\sum_{j=1}^{l-1} \frac{(1/2)_j}{(d/2)_j} \frac{C_{(j)}(\Sigma)}{j!} &= R_1(\Sigma) - R_l(\Sigma) \\
&= O(d^{-(1-r)/2}) + O(d^{-l(1-r)/2}) \\
&= O(d^{-(1-r)/2}),
\end{align*}
and therefore 
\begin{align*}
\bigg(\sum_{j=1}^{l-1} \frac{(1/2)_j}{(d/2)_j} \frac{C_{(j)}(\Sigma)}{j!}&\bigg) \cdot O(d^{-[1+(m-1)(1-r)]/2}) \\
&= O(d^{-(1-r)/2}) \cdot O(d^{-[1+(m-1)(1-r)]/2}) \\
&= O(d^{-[1+m(1-r)]/2}).
\end{align*}
The last $O(\cdot)$ term in \eqref{eq_CovX_BigO} is 
$$
O(d^{-(1-r)}) \cdot O(d^{-[1+(m-1)(1-r)]/2}) = O(d^{-[1+(m+1)(1-r)]/2}).
$$
Collecting all $O(\cdot)$ terms in \eqref{eq_CovX_BigO}, we obtain 
\begin{align*}
\Cov(X) &= \bigg(1 - \sum_{j=1}^{m-1} \frac{(1/2)_j}{(d/2)_j} \frac{C_{(j)}(\Sigma)}{j!}\bigg) \bigg(\sum_{k=1}^{m-1} \frac{(1/2)_k}{(d/2)_k} \frac{\nabla C_{(k)}(\Sigma)}{k!}\bigg) \nonumber \\
&\qquad + O(d^{-(3-2r)/2}) + O(d^{-[1+(m-1)(1-r)]/2}) \\
&\qquad + O(d^{-[1+m(1-r)]/2}) + O(d^{-[1+(m+1)(1-r)]/2}).
\end{align*}
On applying \eqref{eq_sum_of_Big_Os}, we obtain \eqref{eq_cov_matrix_2}.  
$\Box$

\subsection{The proofs of (\ref{eq_ak_bound}) and (\ref{eq_rk_bound_3})}
\label{sec_lemmas}

First, we establish the inequality \eqref{eq_ak_bound} for the coefficients $a_k$ defined in \eqref{eq_ak}.  

\begin{lemma}
\label{lem_ak_formula}
For $k = 0,1,2,\ldots$, 
\begin{equation}
\label{eq_ak_formula}
a_k = \sum_{l=0}^k \frac{\big((d^{1/2}-1)/2\big)^l}{l!} \frac{(1/2)_{k-l}}{(k-l)!} \le \frac{(d^{1/2}/2)_k}{k!}.
\end{equation}
\end{lemma}

\Proof
First, we follow the approach of Comtet \cite[p.~97, Eq.~2d]{Comtet} to derive an explicit formula for $a_k$.  

Let $t = (t_1,t_2,t_3,\ldots)$ be a vector of indeterminates, and define 
$$
a_k(t) = \sum_{i_1 + 2 i_2 + 3 i_3 + \cdots + k i_k = k} \prod_{j=1}^k \frac{t_j^{i_j}}{i_j! \, (2j)^{i_j}},
$$
$k \ge 0$.  For an indeterminate $u$, a formal generating-function for the sequence $\{a_k(t), k=0,1,2,\ldots\}$ is 
\begin{align*}
G_t(u) = \sum_{k=0}^\infty a_k(t) u^k &= \sum_{k=0}^\infty u^k \sum_{i_1 + 2 i_2 + 3 i_3 + \cdots + k i_k = k} \prod_{j=1}^k \frac{t_j^{i_j}}{i_j! \, (2j)^{i_j}} \\
&= \sum_{k=0}^\infty \sum_{i_1 + 2 i_2 + 3 i_3 + \cdots + k i_k = k} \prod_{j=1}^k \frac{t_j^{i_j} u^{ji_j}}{i_j! \, (2j)^{i_j}} \\
&= \prod_{j=1}^\infty \bigg(\sum_{i_j=0}^\infty \frac{t_j^{i_j} u^{ji_j}}{i_j! \, (2j)^{i_j}}\bigg) \\
&= \prod_{j=1}^\infty \exp\Big(\frac{t_j u^j}{2j}\Big) = \exp\bigg(\sum_{j=1}^\infty \frac{t_j u^j}{2j}\bigg).
\end{align*}
Now set $t_1 = d^{1/2}$ and $t_j = 1$ for all $j \ge 2$.  Then $a_k(t)$ reduces to $a_k$; $G_t(u)$ reduces to $G(u)$, the formal generating-function for the sequence $\{a_k\}$; and we obtain 
\begin{align*}
G(u) = \sum_{k=0}^\infty a_k u^k &= \exp\bigg(\tfrac12 d^{1/2} u + \tfrac12 \sum_{j=2}^\infty \frac{u^j}{j}\bigg) \\
&= \exp\bigg(\tfrac12 (d^{1/2}-1) u + \tfrac12 \sum_{j=1}^\infty \frac{u^j}{j}\bigg) \\
&= \exp\big(\tfrac12 (d^{1/2}-1) u - \tfrac12 \log(1-u)\big) \\
&= \exp\big((d^{1/2}-1) u/2\big) \cdot (1-u)^{-1/2}.
\end{align*}
Expanding both of the latter functions in infinite series, we obtain 
\begin{align*}
G(u) &= \bigg(\sum_{l=0}^\infty \frac{\big((d^{1/2}-1)/2\big)^l}{l!} u^l\bigg) \bigg(\sum_{m=0}^\infty \frac{(1/2)_m}{m!} u^m\bigg) \\
&= \sum_{k=0}^\infty u^k \sum_{l=0}^k \frac{\big((d^{1/2}-1)/2\big)^l}{l!} \frac{(1/2)_{k-l}}{(k-l)!}.
\end{align*}
By comparing the coefficients of $u^k$ in $G_t(u)$ and $G(u)$, we obtain the equality in \eqref{eq_ak_formula}.  

Since $\big((d^{1/2}-1)/2\big)^l \le \big((d^{1/2}-1)/2\big)_l$ then, by applying the well-known convolution identity, 
$$
\sum_{l=0}^k \frac{(v_1)_l}{l!} \frac{(v_2)_{k-l}}{(k-l)!} = \frac{(v_1 + v_2)_k}{k!},
$$
$v_1, v_2 \in \R$, we obtain 
$$
a_k \le \sum_{l=0}^k \frac{\big((d^{1/2}-1)/2\big)_l}{l!} \frac{(1/2)_{k-l}}{(k-l)!} 
= \frac{(d^{1/2}/2)_k}{k!}.
$$
The proof of \eqref{eq_ak_formula} now is complete.  
$\Box$

\smallskip

Next, we establish \eqref{eq_rk_bound_3}.

\begin{lemma}
\label{lem_rem_rk_ineq}
Let $\gamma_1 = \tfrac12(\sqrt{3}+1)$.  Then, for all $k \ge 1$ and $d \ge 1$, 
\begin{equation}
\label{eq_rk_bound_1}
\frac{(d^{1/2}/2)_k}{(d/2)_k} \le \gamma_1^{k-1} ((k-1)!)^{1/2} d^{-k/2}.
\end{equation}
\end{lemma}

\Proof
Denote by $r_k$ the left-hand side of \eqref{eq_rk_bound_1}.  Then $r_1 = d^{-1/2}$ and, for $j \ge 1$, 
$$
\frac{r_{j+1}}{r_j} = \frac{(d^{1/2}/2)_{j+1}}{(d/2)_{j+1}} \cdot \frac{(d/2)_j}{(d^{1/2}/2)_j} = \frac{d^{1/2}+2j}{d+2j}.
$$
We claim that 
\begin{equation}
\label{eq_rj_inequality}
\frac{d^{1/2}+2j}{d+2j} \le \gamma_1 j^{1/2} d^{-1/2},
\end{equation}
$j \ge 1$, and we prove this as follows.  Define $y \equiv y(t) = (t^2+2jt)/(t^2+2j)$,  $t \ge 1$.  It is simple to verify that $y(t)$ attains its maximum at $t_j := 1 + (2j+1)^{1/2}$ and that 
$$
y(t_j) = \frac12 \big(1 + (2j+1)^{1/2}\big).
$$
For $j \ge 1$, it is also straightforward to show that 
$$
\frac12 \big(1 + (2j+1)^{1/2}\big) \le \gamma_1 j^{1/2}.
$$
Consequently, 
$$
\frac{d+2jd^{1/2}}{d+2j} = y(d^{1/2}) \le y(t_j) \le \gamma_1 j^{1/2}.
$$
Multiplying the left- and right-hand sides of the above inequality by $d^{-1/2}$, we obtain \eqref{eq_rj_inequality}.  
$\Box$


\end{document}